\documentclass[12pt]{amsart}

\usepackage[textheight=23.5cm, left=1.2cm, right=1.2cm]{geometry}

\usepackage[english]{babel}
\usepackage[T1]{fontenc}
\usepackage[utf8]{inputenc}

\theoremstyle{plain}
    \newtheorem{theorem}{Theorem}[section]
    \newtheorem{lemma}[theorem]{Lemma}
    \newtheorem{proposition}[theorem]{Proposition}
    
    \newtheorem{conjecture}[theorem]{Conjecture}

\theoremstyle{definition}

    \newtheorem{example}[theorem]{Example}

\usepackage[
    draft = false,
    unicode = true,
    colorlinks = true,
    allcolors = blue,
    hyperfootnotes = true,
    citecolor = red
]{hyperref}
\usepackage{amsfonts,amsmath,amssymb,amscd,amsthm,url}
\usepackage[inline]{enumitem}
\usepackage{graphicx,epsf,afterpage}
\usepackage{soul}
\usepackage{multicol}
\usepackage{array}
\usepackage{braket}
\usepackage{epigraph}
\usepackage{rotating}

\usepackage[
backend=biber,
style=numeric-comp,
sorting=none,
giveninits=true
]{biblatex}
\usepackage{csquotes}

\usepackage{xcolor}

\usepackage{ulem}

\usepackage{tikz}
\usetikzlibrary{positioning, calc, intersections, through, arrows, matrix, chains, math}
\usetikzlibrary{shadows}
\usetikzlibrary{graphs}
\usetikzlibrary{graphs.standard}
\usetikzlibrary{patterns}
\usetikzlibrary{decorations.pathreplacing,calligraphy,backgrounds}

\newcommand\norm[1]{\ensuremath{\left\lVert#1\right\rVert}}
\newcommand\abs[1]{\ensuremath{\left\lvert#1\right\rvert}}

\renewcommand{\Pr}{\mathrm{P}}

\DeclareMathOperator{\Expect}{\mathbb{E}}

\DeclareMathOperator{\Img}{Im}

\newcommand{\Lcal}{\mathcal{L}}

\newcommand{\R}{\ensuremath{\mathbb{R}}}

\newcommand{\Cx}{\ensuremath{\mathbb{C}}}

\newcommand{\N}{\ensuremath{\mathbb{N}}}

\renewcommand{\geq}{\geqslant}

\renewcommand{\leq}{\leqslant}

\newcounter{mcnt}

\newcounter{wordcnt}

\begin{document}

\title[The SLLN for random semigroups on uniformly smooth Banach spaces]{The Strong Law of Large Numbers \\ for random semigroups \\ on uniformly smooth Banach spaces}

\author{S.\,V. Dzhenzher and V.\,Zh. Sakbaev}

\thanks{\hspace{-4mm}S.\,V. Dzhenzher: sdjenjer@yandex.ru
\\
V.\,Zh. Sakbaev: fumi2003@mail.ru
\\
V.\,Zh. Sakbaev: Keldysh Institute of Applied Mathematics of Russian Academy of Sciences, Moscow, Russia
\\
Both authors: Moscow Institute of Physics and Technology, Dolgoprudny, Russia
}

\begin{abstract}
     We consider random linear continuous operators $\Omega \to \mathcal{L}(\mathcal{X}, \mathcal{X})$ on a Banach space $\mathcal{X}$.
     For example, such random operators may be random quantum channels.
     The Law of Large Numbers is known when $\mathcal{X}$ is a Hilbert space, in the form of the usual Law of Large Numbers for random operators, and in some other particular cases.
     Instead of the sum of i.i.d. variables, there may be considered the composition of random semigroups $e^{A_i t/n}$.
     We obtain the Strong Law of Large Numbers in Strong Operator Topology for random semigroups of bounded linear operators on a uniformly smooth Banach space.
     We also develop another approach giving the SLLN in Weak Operator Topology for all Banach spaces.
\end{abstract}

\maketitle

\thispagestyle{empty}

\emph{Keywords: Random semigroups; Strong Law of Large Numbers; Uniformly smooth Banach spaces}

\vspace{5mm}

\emph{MSC: 28B05, 47D03, 60F15}

\tableofcontents

\section{Introduction}

The theory of compositions of random linear operators is constructed and developed in a lot of publications \cite{Mehta, Furstenberg, Tutubalin}. 
A random quantum dynamics can be considered either as a composition of random linear operators in a Hilbert space \cite{Kempe, OSS19}, or as a composition of random quantum channels \cite{Aaron, Hol, GOSS-2022, New-23, N-23}.
See in \cite{Berezin, OSS-2016, Orlov2, Aaron, Hol, Pechen-Volovich, GOSS21, GOSS-2022, Joye, Pechen} how random operator-valued processes arise.
See detailed historical exposition in \cite{DzhenzherSakbaev24, DzhenzherSakbaev25}.

The Law of Large Numbers (LLN) for random variables acting in a Banach space is well studied \cite{Ledoux-Talagrand-1991-probability}.
Distributions of products of random linear operators have been studied in \cite{Oseledets, Tutubalin, Skorokhod, New-23}.
The LLN for compositions of random operators in different operator topologies, which is the non-commutative generalization of the LLN for random variables, was studied in \cite{S16, S18, OSS19}. Different versions of limit theorems of LLN and CLT type for compositions of i.i.d. random operators were obtained in \cite{Berger, GOSS-2022, SSSh23}.
The averaging of random quantum channels is studied in \cite{AmosovSakbaev}, where the conditions for the invariance of the set of normal states w.r.t. an average quantum channel are obtained, and the conditions for the destruction of the set of normal states under the average quantum channel are presented.

A lot of results on the LLN were obtained \cite{S16, S18, OSS19} for random linear operators acting in a Hilbert space. But a consideration of quantum dynamics of an open quantum system needs to study the random operators in the Banach space of quantum states \cite{Hol, GOSS-2022, Kalmetev-Orlov-Sakbaev-2023, Pechen, Sahu}. However, the case of random operators acting in a Banach space is studied significantly less. This paper tries to fill this gap for the Strong Law of Large Numbers (SLLN) for random semigroups of bounded linear operators in a Banach space; see some results on the LLN for a Banach space of linear operators acting in a finite dimensional space in \cite{GOSS-2022}, for a Banach space \(\ell_p,\,p\in [\,1,2\,],\) of sequences in \cite{DzhenzherSakbaev24}, and for a $p$-th Schatten class, \(p\in[\,1,2\,]\), with a discrete probability space, in \cite{DzhenzherSakbaev25}.

The structure of the article is the following.
In~\S\ref{s:main-res}, the
main results (Theorems~\ref{t:slln-sot}, \ref{t:slln-wot}) are given.
In~\S\ref{s:lemmas}, the auxiliary lemmas are given.
In~\S\ref{s:proof-slln-sot}, the proof of Theorem~\ref{t:slln-sot} (SLLN in SOT) is given.
In~\S\ref{s:proof-slln-wot}, the proof of Theorem~\ref{t:slln-wot} (SLLN in WOT) is given with the development of another approach giving the SLLN for seminorms induced by sesquilinear forms.

\section{Background and statements}\label{s:main-res}

\subsection*{Topology}

In this paper, we consider Banach spaces over some field $\mathbb{K}$ ($\R$ or $\Cx$).
Recall that for Banach spaces $\mathcal{X},\mathcal{Y}$, the space $\Lcal(\mathcal{X},\mathcal{Y})$ is the Banach space of linear bounded operators $\mathcal{X} \to \mathcal{Y}$.
We shortly denote \(\Lcal(\mathcal{X}) := \Lcal(\mathcal{X},\mathcal{X})\).

For a Banach space $\mathcal{X}$ we denote by $\mathcal{X}^*$ the \emph{dual space} of continuous linear functionals \(\mathcal{X}\to\mathbb{K}\).
For \(x \in \mathcal{X}\) and \(f \in \mathcal{X}^*\) we identify \(\braket{f, x} := f(x)\).
We consider that the complex conjugate is implemented in the action of the dual element; so \(\braket{,}\) is a sesquilinear form.

A Banach space \(\mathcal{X}\) is said to be \textbf{\emph{uniformly smooth}} if the limit
\[
    \lim_{t\to 0} \frac{\norm{x + ty} - \norm{x}}{t}
\]
exists for all \(x,y\) from the unit sphere in \(\mathcal{X}\).
Following \cite{Pisier-2016-martingales} we say that $\mathcal{X}$ is \textbf{\emph{$p$-smooth}} for some \(p \in (1,2\,]\), if there is $C >0$ such that for any \(x,y\in \mathcal{X}\) we have
\[
    \norm{x+y}^p + \norm{x-y}^p \leq 2\norm{x}^p + C\norm{y}^p.
\]
By \cite[Theorem~4.24]{Pisier-2016-martingales}, any uniformly smooth Banach space is $p$-smooth for some \(p \in (1,2\,]\), and vice versa.

It is well known \cite{BCL-Smoothness} that all Hilbert spaces are $2$-smooth, and that for \(1 < p < \infty\) the spaces \(L_p\) are uniformly smooth.
In particular, the spaces \(\ell_p\) of sequences for \(1 < p < \infty\) are uniformly smooth.
It is proved in \cite{BCL-Smoothness, Kazakh-Smoothness} that $p$-th Schatten class operators over separable Hilbert space for \(1 < p < \infty\) are uniformly smooth.
See also \cite{Ledoux-Talagrand-1991-probability} for the connection of the SLLN with the type and the smoothness of a Banach space.

\subsection*{Probability}

Let \((\Omega, \mathcal{F}, \Pr)\) be a probability space with a complete $\sigma$-additive measure $\Pr$.

For a Banach space $\mathcal{X}$ we say that $\xi\colon\Omega\to \mathcal{X}$ is \textbf{\emph{simple}} if its image is finite and the preimage of every point is measurable.

For a Banach space $\mathcal{X}$ we say that $\xi\colon\Omega\to \mathcal{X}$ is a \textbf{\emph{random element}} if it is \emph{Bochner measurable}, i.e. there is a sequence of simple \(\xi_n\) converging to $\xi$ uniformly almost surely (u.a.s.).
By Pettis's measurability theorem \cite[Theorem~II.1.2]{DiestelUhl} this is equivalent to the fact that $\xi$ is
\begin{enumerate}[label=\arabic*.]
    \item \emph{weakly measurable}, i.e. $\braket{f,\xi}$ is measurable for any \(f\in \mathcal{X}^*\), and
    
    \item \emph{essentially separably valued}, i.e. there is \(N \in \mathcal{F}\) such that \(\Pr(N) = 0\) and \(\xi(\Omega\setminus N)\) is a (norm) separable subset of $\mathcal{X}$ (one may say, $\xi$ almost surely lies in some separable subset of $\mathcal{X}$).
\end{enumerate}
This is well-known \cite{DiestelUhl, Yosida1995} that the norm of a random element is measurable.

In the case when $\mathcal{X}$ is the space of linear bounded operators, acting in some Banach spaces, we will call random elements in $\mathcal{X}$ as \textbf{\emph{random operators}}.

\begin{lemma}[Composition]\label{l:compos}
    Let \(A_n\colon\Omega\to\Lcal(\mathcal{X}, \mathcal{Y})\) and \(B_n\colon\Omega\to\Lcal(\mathcal{Y}, \mathcal{Z})\) be random operators converging u.a.s. to $A$ and $B$, correspondingly.
    Then \(B_nA_n \xrightarrow[n\to\infty]{u.a.s.} BA\).

    Consequently, a composition of random operators is a random operator.
\end{lemma}
\begin{proof}
    We have
    \[
        \norm{B_nA_n - BA} \leqslant \norm{B_nA_n - BA_n} + \norm{BA_n - BA} \leq \norm{B_n-B}\norm{A_n} + \norm{B}\norm{A_n-A} \xrightarrow[n\to\infty]{} 0.
    \]
\end{proof}

A sequence $\{\xi_n\}$ of random elements \emph{converges in $L_1$} to a random $\xi$, denoted \(\xi_n\xrightarrow[n\to\infty]{L_1}\xi,\) if
\[
    \Expect \norm{\xi_n-\xi} \xrightarrow[n\to\infty]{} 0.
\]
A random element $\xi$ is \emph{Bochner integrable} if there is a sequence of simple random elements converging to $\xi$ both u.a.s. and in $L_1$.

For a simple $\xi\colon\Omega\to \mathcal{X}$ its \textbf{\emph{expected value}} is
\[
    \Expect\xi := \sum_{x \in \Img\xi} x \Pr\{\xi=x\}.
\]
For a Bochner integrable $\xi\colon\Omega\to \mathcal{X}$ its \textbf{\emph{expected value}}
\[
    \Expect\xi := \lim_{n\to\infty} \Expect \xi_n
\]
is the Bochner integral, where $\xi_n$ are simple converging to $\xi$ both u.a.s. and in $L_1$.
It is well-known \cite[Theorem~II.2.2]{DiestelUhl} that $\xi$ is Bochner integrable if and only if \(\norm{\xi}\) is integrable.

\begin{lemma}[Integration]
\label{l:integ}
    Let \(A \colon \Omega \to \Lcal(\mathcal{X},\mathcal{Y})\) be a (Bochner) integrable random operator.
    Then for any \(x\in \mathcal{X}\)
    \[
        (\Expect A)x = \Expect (Ax).
    \]
\end{lemma}
\begin{proof}
    It is obvious that $Ax$ is a random element. It is integrable since \(\norm{Ax}\) is.
    
    For a simple $A$ we have
    \[
        (\Expect A)x = \sum_{U \in \Img A} Ux\Pr\{A=U\} = \sum_y y\Pr\{Ax=y\} = \Expect (Ax).
    \]

    In the general case, taking \(A_n \xrightarrow{L_1,\,u.a.s.} A\) we have
    \[
        (\Expect A)x = \lim_{n\to\infty} (\Expect A_n)x = \lim_{n\to\infty} \Expect (A_nx) = \Expect (Ax),
    \]
    where the last equality follows since
    \(
        \Expect \norm{A_nx - Ax} \leq \Expect \norm{A_n - A}\norm{x} \xrightarrow[n\to\infty]{} 0.
    \)
\end{proof}

A function \(\xi\colon\Omega\to \mathcal{X}\) is \textbf{\emph{ball measurable}} if the preimage of any ball is measurable, i.e. \(\{\xi \in B\} \in \mathcal{F}\) for any ball $B \subset \mathcal{X}$.
A function \(\xi\colon\Omega\to \mathcal{X}\) is \textbf{\emph{Borel measurable}} if the preimage of any open set is measurable, i.e. \(\{\xi \in O\} \in \mathcal{F}\) for any open $O \subset \mathcal{X}$.
It is obvious that for essentially separably valued functions these definitions are equivalent.

The following lemma is known \cite{Varadarajan}.

\begin{lemma}[Measurability]\label{l:measure}
    An essentially separably valued function \(\xi\colon\Omega\to \mathcal{X}\) is ball measurable if and only if it is a random element.
\end{lemma}
\begin{proof}
    Let us prove the <<if>> part.
    Take a sequence of simple \(\xi_n \xrightarrow[n\to\infty]{u.a.s.} \xi\). Take \(B = \{y \in \mathcal{X} : \norm{y-x} \leqslant r\}\) for some $x\in \mathcal{X}$ and \(r>0\). Then
    \[
        \{\xi \in B\} = \bigcap_{\varepsilon\to0} \bigcup_{N=1}^\infty \bigcap_{n=N}^\infty \{\norm{\xi_n - x} < r+\varepsilon\} \in \mathcal{F},
    \]
    where the first equality is up to a set of zero probability.

    Let us prove the <<only if>> part.
    Take a separable subset $\mathcal{X}_0 \subset \mathcal{X}$ in which $\xi$ lies almost surely.
    For each $n \in \N$ cover $\mathcal{X}_0$ with a countable family of balls \(B_{n,k}\) of radius $\frac{1}{n}$.
    For each $n$ modify these \(B_{n,k}\) so that \(B_{n,k+1}\) is obtained by deleting \(B_{n,k}\). They become a disjoint cover of $\mathcal{X}_0$.
    Construct countably valued functions $\xi_n$ by setting on \(\xi^{-1}(B_{n,k})\) some points from \(B_{n,k}\), and on the remaining set of zero probability some fixed point.
    Then by \cite[Corollary~II.1.3]{DiestelUhl} \(\xi_n \xrightarrow[n\to\infty]{u.a.s.}\xi\) implies Bochner measurability of $\xi$.     
\end{proof}

Random elements \(\{\xi_n\colon\Omega\to \mathcal{X}_n\}\) are \textbf{\emph{independent}} if for any $n\in \N$ and any open \(O_k \subset \mathcal{X}_k\) we have
\[
    \Pr\{\xi_1 \in O_1,\,\ldots,\,\xi_n \in O_n\} = \Pr\{\xi_1 \in O_1\}\ldots\Pr\{\xi_n\in O_n\}.
\]
Note that one may take $O_k=\mathcal{X}_k$ and so make some events be the whole $\Omega$.

\subsection*{Semigroups}

For a bounded random operator \(A\colon\Omega\to\Lcal(\mathcal{X})\) denote by \(e^{At}\) the random (uniformly continuous) semigroup with $A$ as the infinitesimal generator.

\begin{theorem}[SLLN in SOT; proved in~\S\ref{s:proof-slln-sot}]\label{t:slln-sot}
    Let \(\mathcal{X}\) be a uniformly smooth Banach space.
    Let \(A\colon\Omega\to\Lcal(\mathcal{X})\) be a generator of a semigroup $e^{At}$ such that $\norm{A}$ is bounded.
    Let $A,A_1, A_2, \ldots$ be a sequence of i.i.d. random generators.
    Then $e^{A_1t/n} \ldots e^{A_nt/n}$ converges a.s. in SOT of $\norm{\cdot}_\mathcal{X}$ to $e^{\Expect At}$ uniformly for $t>0$ in any segment; that is, for any \(x \in \mathcal{X}\) and $T > 0$ almost surely
    \[
        \lim_{n\to\infty} \sup\limits_{t \in [\,0,T\,]} \norm{(e^{A_1t/n} \ldots e^{A_nt/n} - e^{\Expect At})x} = 0.
    \]
\end{theorem}

There is an example \cite[Example~5]{OSS19} when the LLN fails even for a Hilbert space \(\mathcal{X}\) in the case of a finitely additive measure (but not $\sigma$-additive).
In that case, Bochner measurability is replaced with weak measurability, and expected value becomes the Pettis integral.

It is interesting to notice that in the usual Strong Law of Large Numbers for random variables it is necessary to demand only that all variables are uniformly bounded by some constant, independent and have the same expected value.
So in the classical case the requirement of the random variables being i.i.d. is redundant.
However, it is not clear if this requirement is redundant in the case of random semigroups of operators.

\begin{theorem}[SLLN in WOT; proved in~\S\ref{s:proof-slln-wot}]\label{t:slln-wot}
    Let $\mathcal{X}$ be a Banach space.
    Then under the assumptions of Theorem~\ref{t:slln-sot}, $e^{A_1t/n} \ldots e^{A_nt/n}$ converges a.s. in WOT to $e^{\Expect At}$ uniformly for $t>0$ in any segment; that is, for any \(x \in \mathcal{X}\), \(f \in \mathcal{X}^*\) and $T > 0$ almost surely
    \[
        \lim_{n\to\infty} \sup\limits_{t \in [\,0,T\,]} \abs{\braket{f,(e^{A_1t/n} \ldots e^{A_nt/n} - e^{\Expect At})x}} = 0.
    \]
\end{theorem}

Thus by the Schur's Theorem we immediately obtain that for \(\mathcal{X}=\ell_1\) the SLLN in SOT of \(\norm{\cdot}_{\ell_1}\) holds.
Note that \(\ell_1\) is \emph{not} uniformly smooth.
This leads us to the following conjecture.

\begin{conjecture}
    Let \(\mathcal{X}\) be a Banach space.
    Then under the assumptions of Theorem~\ref{t:slln-sot}, $e^{A_1t/n} \ldots e^{A_nt/n}$ converges a.s. in SOT of $\norm{\cdot}_\mathcal{X}$ to $e^{\Expect At}$ uniformly for $t>0$ in any segment.
\end{conjecture}

\section{Auxiliary lemmas and notions}\label{s:lemmas}

\begin{lemma}[Independence]\label{l:indep}
    Let \(A_1, \ldots, A_k\) be independent random operators with bounded norms such that their composition is well-defined (so that their images and domains are correspondent).
    Then
    \[
        \Expect (A_1\ldots A_k) = (\Expect A_1)\ldots (\Expect A_k).
    \]
\end{lemma}
\begin{proof}
    If all \(A_1, \ldots, A_k\) are simple, then
    \[
        \Expect (A_1\ldots A_k) =
        \sum_{U_1, \ldots, U_k} U_1\ldots U_k \Pr\{A_1=U_1\}\ldots\Pr\{A_k=U_k\} =
        (\Expect A_1)\ldots (\Expect A_k).
    \]
    In the general case, take \(A_{n,i} \xrightarrow[n\to\infty]{L_1,\,u.a.s.}A_i\) such that $\{A_{n,i}\}_{i=1,\ldots,k}$ are independent for each $n \in \N$. This can be done since $\sigma$-algebras \(\sigma(A_i)\) are independent, and since $A_{n,i}$ can be constructed ball-measurable w.r.t. corresponding $\sigma$-algebras \(\sigma(A_i)\).
    Then
    \[
        \Expect (A_1\ldots A_k) = \lim_{n\to\infty} \Expect(A_{n,1}\ldots A_{n,k}) =
        \lim_{n\to\infty} \Expect(A_{n,1})\ldots \Expect(A_{n,k}) =
        (\Expect A_1)\ldots (\Expect A_k).
    \]
\end{proof}

The ideas of usage and proofs of the following two lemmas are due to \cite{GOSS-2022}.

\begin{lemma}[Chernoff's equivalence \cite{Chernoff-1968}]\label{l:chernoff-eq}
    Let \(\mathcal{X}\) be a Banach space.
    Let \(A\colon\Omega\to\Lcal(\mathcal{X})\) be a generator of a semigroup $e^{At}$ such that $\norm{A}$ is bounded.
    Let $A,A_1, A_2, \ldots$ be a sequence of i.i.d. random generators.
    Then for any \(x \in \mathcal{X}\) and $T > 0$
    \[
        \lim\limits_{n\to\infty} \sup\limits_{t \in [\,0,T\,]} \norm{\left(e^{\Expect At} - \Expect e^{A_1t/n} \ldots e^{A_nt/n}\right)x} = 0.
    \]
\end{lemma}

\begin{proof}
    Denote the~mapping $F\colon [\,0, +\infty)\to\Lcal(\mathcal{X})$ by the rule $F(t) := \Expect e^{At}$.
    By the Chernoff's theorem we have
    \[
        \lim\limits_{n\to\infty} \sup\limits_{t \in [\,0,T\,]} \norm{\left(e^{F'(0)t} - F\left(t/n\right)^n\right)x} = 0
    \]
    if (a) $F$ is continuous,
    (b) \(F(0) = I\),
    (c) \(\norm{F(t)} \leqslant e^{at}\) for some $a > 0$, and
    (d) for any $y\in X$ there is \(F'(0)y := \lim\limits_{t \to +0} \dfrac{F(t)y - y}{t}\).
    Here~(a,b) are obvious,
    (c) holds since \(\norm{A}\) is bounded,
    and~(d) holds for $F'(0) = \Expect A$.
    It remains to notice that the Chernoff's theorem gives exactly the required equality since
    \(F'(0) = \Expect A\) and \(\Expect e^{A_1t/n} \ldots e^{A_nt/n} = F(t/n)^n\).
\end{proof}

Let \(\mathcal{X}\) be a Banach space and $A$ be an integrable operator.
Let \(A,A_1,\ldots,A_n\) be integrable operators with bounded norms.
For $i \in \N$ and $s \geq 0$ denote
\[
    \Delta_i(s) := e^{A_is} - \Expect e^{As}.
\]
For integer $0 \leq k \leq n$ denote \([n]:= \{1,\ldots,n\}\), and denote by \(\binom{[n]}{k}\) the family of $k$-element subsets of \([n]\). 
For \(\{i_1 < \ldots < i_k\} \in \binom{[n]}{k}\) and \(s>0\) denote
\[
    F_{n,\{i_1, \ldots, i_k\}}(s) := (\Expect e^{As})^{i_1-1} \Delta_{i_1}(s) (\Expect e^{As})^{i_2-i_1-1} \ldots \Delta_{i_k}(s) (\Expect e^{As})^{n-i_k},
\]
meaning
\[
    F_{n,\varnothing}(s) = (\Expect e^{As})^n.
\]

\begin{lemma}\label{l:est-norm}
    Let \(\mathcal{X}\) be a Banach space.
    Let $A, A_1, A_2, \ldots, A_n$ be i.i.d. random operators with their norms being bounded by some \(\rho > 0\).
    Then for each $n\in \N$ and \(P \in \binom{[n]}{k}\)
    \[
        \norm{F_{n,P}(s)} \leq (2\rho s)^k e^{n\rho s}.
    \]
\end{lemma}

\begin{proof}
    We have \(\Expect\norm{e^{As}} \leq e^{\rho s}\).
    Hence by \cite[Theorem~II.4.ii]{DiestelUhl}
    \begin{equation}\label{eq:e-norm}
        \norm{\Expect e^{As}} \leq e^{\rho s}.
    \end{equation}
    Denote by \(I\colon\mathcal{X}\to\mathcal{X}\) the identity operator.
    By \cite[Corollary~4.6.4 (the mean value theorem)]{BogachevSmolyanov}
    \[
        \norm{e^{As} - I} \leq \rho se^{\rho s}.
    \]
    By~\eqref{eq:e-norm} and again by \cite[Corollary~4.6.4 (the mean value theorem)]{BogachevSmolyanov}
    \[
        \norm{\Expect e^{As} - I} \leq \rho se^{\rho s}.
    \]
    Hence by the triangle inequality
    \begin{equation}\label{eq:delta-norm}
        \norm{\Delta_i(s)} \leq 2\rho s e^{\rho s}.
    \end{equation}
    Now the lemma follows as the binomial combinations of inequalities~(\ref{eq:e-norm}, \ref{eq:delta-norm}).
\end{proof}

\section{Proof of Theorem~\ref{t:slln-sot} (SLLN in SOT)}\label{s:proof-slln-sot}

For an integrable random element \(\xi\) and a sub-$\sigma$-algebra \(\mathcal{G}\subset\mathcal{F}\) one may define the \emph{conditional expectation} \(\Expect^\mathcal{G}\xi\) analogously to the real-valued case; see \cite{DiestelUhl} for details.

The sequence \((\mu)_{n\in\N}\) of random elements is called \emph{a martingale relative to a filtration \((\mathcal{F}_{n\in\N})\)} if each \(\mu_n\) is \(\mathcal{F}_n\)-measurable and \(\Expect^{\mathcal{F}_{n-1}}\mu_n = \mu_{n-1}\); see \cite{DiestelUhl, Pisier-2016-martingales} for details.

The crucial idea of the proof and of the usage of the uniform smoothness is the following Burkholder-type theorem \cite[Theorem~4.52]{Pisier-2016-martingales}.

\begin{theorem}[Burkholder-type]\label{t:burkholder}
    Let $\mathcal{X}$ be a $p$-smooth Banach space.
    Let \(\mu_n = \sum_{k=1}^n d_k\) be a martingale.
    Then for any \(r \geq 1\) there is some \(C = C(\mathcal{X},r,p)>0\) such that
    \[
        \Expect\norm{\mu_n}^r \leq C\Expect\left(\sum_{k=1}^n \norm{d_k}^p\right)^{\frac{r}{p}}.
    \]
\end{theorem}

\noindent
See \cite{Pinelis-1994, Pinelis-2015} for the optimal bounds in Burkholder-type and Rosenthal-type inequalities in $2$-smooth Banach spaces.

The general idea of the proof of Theorem~\ref{t:slln-sot} is due to \cite[proof of Theorem~2]{GOSS-2022}.
There firstly implicitly martingales appear, though were not used.

\begin{proof}[Proof of Theorem~\ref{t:slln-sot} (SLLN in SOT)]
    For a fixed \(x\in\mathcal{X}\) denote
    \[
        \mu_n(t) := \left(e^{A_1t/n} \ldots e^{A_nt/n} - \Expect e^{A_1t/n} \ldots e^{A_nt/n}\right)x =
        \sum_{k=1}^n\sum_{P\in\binom{[n]}{k}} F_{n,P}(t/n)x.
    \]
    By Chernoff's equivalence~\ref{l:chernoff-eq} it is sufficient to prove that for any \(T>0\) a.s.
    \[
        \lim\limits_{n\to\infty} \sup\limits_{t \in [\,0,T\,]} \norm{\mu_n(t)} = 0.
    \]
    For \(k \leq n\) denote
    \[
        d_{n,k}(t) := \sum_{\substack{P \subset [n] \\ \max P = k}} F_{n,P}(t/n)x.
    \]
    Up to \(\norm{x}\) by Lemma~\ref{l:est-norm} we have
    \begin{equation}\label{eq:mart-dif}
        \norm{d_{n,k}(t)} \leq \sum_{j=1}^k \binom{k-1}{j-1} \left(\frac{2\rho t}{n}\right)^j e^{\rho t} =
        \frac{2\rho t}{n} e^{\rho t} \left(1+\frac{2\rho t}{n}\right)^{k-1} \leq \frac{2\rho t}{n}e^{3\rho t}.
    \end{equation}
    Note that for fixed $t>0$ and \(n\in\N\) the sequence of
    \[
        \mu_m(t) = \sum_{k=1}^m d_{n,k}(t),
        \quad\text{where}\quad m\leq n,
    \]
    is a martingale w.r.t. \(\mathcal{F}_m = \sigma(A_1,\ldots,A_m)\), since
    \[
        \Expect^{\mathcal{F}_{m-1}} F_{n,\{i_1 < \ldots < i_j=m\}}(s) =
        (\Expect e^{As})^{i_1-1} \Delta_{i_1}(s) (\Expect e^{As})^{i_2-i_1-1} \ldots \Delta_{i_{j-1}}(s) (\Expect e^{As})^{i_j-i_{j-1}-1} \underbrace{\Expect \Delta_{i_j}(s)}_{=0} (\Expect e^{As})^{n-i_j} = 0.
    \]
    Then by Burkholder-type Theorem~\ref{t:burkholder} for \(r=\frac{2p}{p-1}\) and by inequality~\eqref{eq:mart-dif} up to the unimportant constant we have
    \[
        \Expect\norm{\mu_n(t)}^r \leq \left(n \left(\frac{2\rho t}{n}\right)^pe^{3p\rho t}\right)^{\frac{r}{p}} = \frac{1}{n^2} (2\rho t)^r e^{3r\rho t}.
    \]
    Thus by the Markov's inequality (up to the constant from Theorem~\ref{t:burkholder})
    \[
        \Pr\left\{ \norm{\mu_n(t)} > \varepsilon \right\} \leq \frac{\norm{x}^r}{\varepsilon^r} \frac{1}{n^2} (2\rho t)^r e^{3r\rho t}.
    \]
    Finally by the continuity from below of the probability we obtain
    \[
        \Pr\left\{\sup\limits_{t \in [\,0,T\,]} \norm{\mu_n(t)} > \varepsilon \right\} = O\left(\frac{1}{n^2}\right).
    \]
    Now the desired result follows by the Borel--Cantelli lemma.
\end{proof}

Note that in the application of Theorem~\ref{t:burkholder} it is sufficient to take any \(r > \frac{p}{p-1}\) to complete the proof.

\section{Proof of Theorem~\ref{t:slln-wot} (SLLN in WOT)}\label{s:proof-slln-wot}

In this section we develop an approach giving the SLLN in WOT.

Let $\mathcal{X}$ be a Banach space equipped with a \emph{positive} operator \(i\in \Lcal(\mathcal{X}, \mathcal{X}^*)\), i.e. such that \(\braket{i(x), x} \geqslant 0\) for all $x \in \mathcal{X}$.
For \(x \in \mathcal{X}\) denote
\[
    \norm{x}_i := \sqrt{\braket{i(x),x}}.
\]

\begin{proposition}
    Let $\mathcal{X}$ be a Banach space equipped with a positive operator \(i\in \Lcal(\mathcal{X}, \mathcal{X}^*)\).
    Then \(\norm{\cdot}_i\) is the seminorm on $\mathcal{X}$.
\end{proposition}

\begin{proof}
    Absolute homogeneity follows since \(\braket{,}\) is the sesquilinear form.
    Subadditivity follows since for any \(x,y \in \mathcal{X}\) we have
    \begin{multline*}
        \norm{x+y}_i^2 = \braket{i(x),x}+\braket{i(x),y}+\braket{i(y),x}+\braket{i(y),y} =
        \braket{i(x),x}+2\Re\braket{i(x),y}+\braket{i(y),y} \leq \\ \leq
        \braket{i(x),x}+2\sqrt{\braket{i(x),x}\braket{i(y),y}}+\braket{i(y),y} = \left(\norm{x}_i + \norm{y}_i\right)^2,
    \end{multline*}
    where the inequality is the Cauchy--Schwarz inequality, and the second equality holds since the sesquilinear form \(\braket{i(x),y}\) is Hermitian by the positivity of $i$.
\end{proof}

\begin{example}
    If $\mathcal{X}$ is isomorphic to a Hilbert space, then one may take $i$ given by the Riesz representation theorem and obtain \(\norm{\cdot}_i\) as just the Hilbert space norm.
    Other examples:
    \begin{itemize}
        \item \(\mathcal{X} = \ell_p\) for \(p \in [\,1,2\,]\) is the space of sequences, and $i$ is the inclusion \(\ell_p \hookrightarrow \ell_{\frac{p}{p-1}}\);
    
        \item \(\mathcal{X} = L_p(\mu)\) for \(p \in [\,2,+\infty)\) is the space of $\mu$-measurable functions for which the $p$-th power of the absolute value is Lebesgue integrable over finite $\mu$, and $i$ is the inclusion \(L_p(\mu) \hookrightarrow L_{\frac{p}{p-1}}(\mu)\);
    
        \item \(\mathcal{X} = \mathcal{T}_p\) for \(p \in [\,1,2\,]\) is the Schatten-class space of operators on a separable Hilbert space, and $i$ is the inclusion \(\mathcal{T}_p \hookrightarrow \mathcal{T}_{\frac{p}{p-1}}\).
    \end{itemize}
    In all these cases we have \(\norm{\cdot}_i = \norm{\cdot}_2\).
\end{example}

Note that we always have \(\norm{x}_i \leq \sqrt{\norm{i}}\norm{x}\); it is interesting to find out when \(\norm{\cdot}_i\) and \(\norm{\cdot}_\mathcal{X}\) become equivalent.
The seminorm \(\norm{\cdot}_i\) becomes the norm if and only if $i$ is positive definite; in that case \((\mathcal{X},\norm{\cdot}_i)\) becomes a Hilbert space.

\begin{example}
    Take the Banach space \(\mathcal{X}=\R^2\) with the maximum norm.
    It is not a Hilbert space, but the operator \(i\colon x\mapsto (y \mapsto x_1y_1+x_2y_2)\) induces the Hilbert space \((\R^2, \norm{\cdot}_i)\) with the Euclidean norm.
\end{example}

\begin{example}
    For any Banach space $\mathcal{X}$ with a
    Schauder basis, there are infinitely many non-trivial positive operators \(i \in \Lcal(\mathcal{X}, \mathcal{X}^*)\).
    Indeed, w.l.o.g. assume that the Schauder basis \(\{e_n\} \subset \mathcal{X}\) is normalized.
    Denote by \(\{e_n^*\} \subset \mathcal{X}^*\) the biorthogonal functionals.
    Fix any integer \(N \geqslant 1\) and define
    \[
        i\left(\sum_{n=1}^\infty \alpha_ne_n \right) := \sum_{n=1}^N   {\alpha_n}e_n^*.
    \]
    Then operator $i$ is positive since \(\braket{i(x), x} = \sum_{n=1}^N \abs{\alpha_n}^2 \geqslant 0\) for any $x = \sum_n \alpha_n e_n \in \mathcal{X}$.
\end{example}



The following theorem is known in the weaker form of convergence in probability in the case of unitary semigroups generated by self-adjoint operators.
We consider the case of bounded generators; see \cite{DzhenzherSakbaev24, DzhenzherSakbaev25} for simpler cases.

\begin{theorem}\label{t:slln-inorm}
    Let \((\mathcal{X},\norm{\cdot}_i)\) be a seminormed Banach space.
    Let \(A\colon\Omega\to\Lcal(\mathcal{X})\) be a generator of a semigroup $e^{At}$ such that $\norm{A}$ is bounded.
    Let $A,A_1, A_2, \ldots$ be a sequence of i.i.d. random generators.
    Then $e^{A_1t/n} \ldots e^{A_nt/n}$ converges a.s. in SOT of $\norm{\cdot}_i$ to $e^{\Expect At}$ uniformly for $t>0$ in any segment; that is, for any \(x \in \mathcal{X}\) and $T > 0$ almost surely
    \[
        \lim_{n\to\infty} \sup\limits_{t \in [\,0,T\,]} \norm{(e^{A_1t/n} \ldots e^{A_nt/n} - e^{\Expect At})x}_i = 0.
    \]
\end{theorem}

\begin{proof}[Deduction of Theorem~\ref{t:slln-wot} (SLLN in WOT) from Theorem~\ref{t:slln-inorm}]
    Take \(i\in\Lcal(\mathcal{X}, \mathcal{X}^*)\) defined by \(i(y) := f(y)f\). Then \(\norm{y}_i = \abs{f(y)}\), and the result follows by Theorem~\ref{t:slln-inorm}.
\end{proof}

Recall that for a random operator \(A\colon \Omega\to \Lcal(\mathcal{X}, \mathcal{Y})\) the \emph{adjoint random operator} \(A^*\colon \Omega \to \Lcal(\mathcal{Y}^*, \mathcal{X}^*)\) is defined by
\[
    \braket{A^*y, x} = \braket{y, Ax}
\]
for any \(x\in \mathcal{X}\), \(y\in \mathcal{Y}^*\).
The following lemma provides the correctness of the definition.

\begin{lemma}[Adjoint]\label{l:adjoint}
    Let \(A\colon\Omega\to\Lcal(\mathcal{X},\mathcal{Y})\) be a random operator.
    Then its adjoint $A^*$ is also a random operator.
    If additionally $A$ is (Bochner) integrable then so is $A^*$, and
    \[
        \Expect A^* = (\Expect A)^*.
    \]
\end{lemma}
\begin{proof}
    The adjoint $A^*$ is Bochner measurable since for simple \(A_n \xrightarrow[n\to\infty]{u.a.s.} A\) we have \(A_n^* \xrightarrow[n\to\infty]{u.a.s.} A^*\), and $A_n^*$ are also simple for all \(n\in\N\).
    It is Bochner integrable since \(\norm{A^*}=\norm{A}\).
    Take any \(x\in \mathcal{X}\) and \(f\in \mathcal{Y}^*\).
    Consider the linear functional \(T_x\in \mathcal{Y}^{**}\) defined by \(T_xf := \braket{f,x}^*\), where the last star is the complex conjugate.
    Since \(A^*f\) is Bochner integrable, so is \(T_x (A^*f)\), and
    \[
        T_x \Expect (A^*f) = \Expect T_x(A^*f).
    \]
    Then
    \[
        \braket{(\Expect A^*)f, x} = \Expect\braket{A^*f,x} = \Expect\braket{f, Ax} = \braket{f, (\Expect A)x} = \braket{(\Expect A)^*f,x},
    \]
    where the equality before the last follows analogously (or since the expected value is the Pettis integral).
\end{proof}

\begin{lemma}[Random element]\label{l:rand-el}
    Let \(\xi\colon\Omega\to \mathcal{X}\) and \(A\colon\Omega\to\Lcal(\mathcal{X},\mathcal{Y})\) be independent random element and random operator with bounded norms.
    Then
    \[
        \Expect (A\xi) = (\Expect A)(\Expect\xi).
    \]
\end{lemma}
\begin{proof}
    The measurability of \(A\xi\) can be proved analogously to Lemma~\ref{l:compos} (Composition); it has the bounded norm and hence (Bochner) integrable.

    For simple \(A,\xi\) we have
    \[
       \Expect (A\xi) = \sum_{U,x} Ux\, \Pr\{A=U,\,\xi=x\} = (\Expect A)(\Expect \xi).
    \]
    In the general case, take simple \(A_n\) and \(\xi_n\) converging u.a.s. to \(A\) and \(\xi\).
    Then \(A_n\xi_n\) converge u.a.s. to \(A\xi\).
    Hence
    \[
        \Expect (A\xi) = \lim_{n\to\infty} \Expect A_n\xi_n =
        \lim_{n\to\infty} (\Expect A_n) (\Expect \xi_n) = (\Expect A) (\Expect \xi).
    \]
\end{proof}

Note that the proof of the last lemma is similar to the proof of Lemma~\ref{l:indep} (Independence) but not the same, since in Lemma~\ref{l:rand-el} we have the action, while in Lemma~\ref{l:indep} we have the composition.

We will use notations of previous sections.

\begin{proof}[Proof of Theorem~\ref{t:slln-inorm}]
    As in the proof of Theorem~\ref{t:slln-sot}, denote
    \[
        \mu_n(t) := \left(e^{A_1t/n} \ldots e^{A_nt/n} - \Expect e^{A_1t/n} \ldots e^{A_nt/n}\right)x =
        \sum_{k=1}^n\sum_{P\in\binom{[n]}{k}} F_{n,P}(t/n)x.
    \]
    Hence
    \begin{multline*}
        \Expect \norm{\mu_n(t)}_i^4 =
        \Expect \Braket{\sum_{k=1}^n\sum_{m=1}^n \sum_{P\in\binom{[n]}{k}} \sum_{Q\in\binom{[n]}{m}} F_{n,P}(t/n)^*iF_{n,Q}(t/n)x,x }^2 = \\ =
        \sum_{k=1}^n\sum_{m=1}^n\sum_{k'=1}^n\sum_{m'=1}^n \sum_{P\in\binom{[n]}{k}} \sum_{Q\in\binom{[n]}{m}} \sum_{P'\in\binom{[n]}{k'}} \sum_{Q'\in\binom{[n]}{m'}}
        \Expect \braket{F_{n,P}(t/n)^*iF_{n,Q}(t/n)x,x} \braket{F_{n,P'}(t/n)^*iF_{n,Q'}(t/n)x,x}.
    \end{multline*}
    In the last sum there survive the summands for which each element of each of the sets $P,Q,P',Q'$ is contained in at least one of three other sets.
    Indeed, suppose w.l.o.g. that some \(j\in P\) does not lie in \(P'\cup Q\cup Q'\).
    Define the super-operator \(V\) that takes an operator $B$ as an argument and maps it to \(\braket{F_{n,P}(t/n)^*iF_{n,Q}(t/n)x,x} \braket{F_{n,P'}(t/n)^*iF_{n,Q'}(t/n)x,x}\), where instead of \(A_j\) the operator $B$ is substituted.
    Then by Lemma~\ref{l:rand-el} (Random element) we have that
    \[
        \Expect (VA_j) = (\Expect V)(\Expect A_j) = 0.
    \]

    For a finite set $P$ denote by \(\abs{P}\) the number of elements in $P$.
    Let $\rho$ be the radius of the ball which bounds the image of $A$, i.e. $\norm{A} \leq \rho$.
    Then by Lemma~\ref{l:est-norm} up to the multiplicative constant \(e^{4\rho t}\norm{x}^4\norm{i}^2\) we get
    \[
        \Expect \norm{\mu_n(t)}_i^4 \leq  \sum_{P,P',Q,Q'} \left(\frac{2\rho t}{n}\right)^{\abs{P}+\abs{P'}+\abs{Q}+\abs{Q'}},
    \]
    where the sum is over tuples \((P,P',Q,Q')\) of subsets of \(\{1,\ldots,n\}\) with the property that each element of each of the sets is contained in the union of the three other sets.
    For \(j\in\{1,\ldots,n\}\) and a tuple \((P,P',Q,Q')\) denote by \(\ell_j(P,P',Q,Q')\in\{0,2,3,4\}\) the number of sets containing the element $j$.
    Then the last sum by the inclusion-exclusion principle equals
    \begin{multline*}
        \sum_{P,P',Q,Q'} \left(\frac{2\rho t}{n}\right)^{\sum_{j=1}^n \ell_j(P,P',Q,Q')} =
        \sum_{P,P',Q,Q'} \prod_{j=1}^n \left(\frac{2\rho t}{n}\right)^{\ell_j(P,P',Q,Q')} = \\ =
        \binom{4}{0}S_4^n - \binom{4}{1}S_3^n + \binom{4}{2}S_2^n - \binom{4}{3}S_1^n + \binom{4}{4}S_0^n,
    \end{multline*}
    where
    \begin{align*}
        S_4 &= \binom{4}{0} + \binom{4}{2}\left(\frac{2\rho t}{n}\right)^2 + \binom{4}{3}\left(\frac{2\rho t}{n}\right)^3 + \binom{4}{4}\left(\frac{2\rho t}{n}\right)^4, 
        &&S_3 = \binom{3}{0} + \binom{3}{2}\left(\frac{2\rho t}{n}\right)^2 + \binom{3}{3}\left(\frac{2\rho t}{n}\right)^3, \\
        S_2 &= \binom{2}{0} + \binom{2}{2}\left(\frac{2\rho t}{n}\right)^2, 
        &&S_1 = S_0 = 1.
    \end{align*}
    Finally, by the Markov's inequality and according to Taylor's theorem in the Peano form of the remainder for the first order, we get
    \[
        \Pr \left\{ \norm{\mu_n(t)}_i > \varepsilon\right\} \leq \frac{e^{4\rho t}\norm{x}^4\norm{i}^2}{\varepsilon^4}
        \left(\frac{12\rho^2t^2}{n^2} + o\left(\frac{1}{n^2}\right)\right).
    \]
    Now the SLLN follows by continuity from below of the probability and by the Borel--Cantelli lemma.
\end{proof}

\printbibliography

\end{document}